\documentclass[10pt,a4paper]{amsart}
\usepackage{graphicx}
\usepackage{amssymb}
\usepackage{amsmath}
\usepackage{amsthm}
\usepackage{amscd}
\usepackage[all,2cell,dvips]{xy}
\UseAllTwocells \SilentMatrices
\newtheorem{thm}{Theorem}[section]

\newtheorem{lem}[thm]{Lemma}

\theoremstyle{definition}

\theoremstyle{remark}
\newtheorem{rem}[thm]{\bf Remark}
\numberwithin{equation}{section}

\begin{document}
\title[A short proof of  HRS-tilting]
{A short proof of  HRS-tilting}
\author[  Xiao-Wu Chen
] {Xiao-Wu Chen}
%\thanks{$^*$ The corresponding author}
\thanks{This project was supported by Alexander von Humboldt Stiftung, and was
also partially supported by China Postdoctoral Science Foundation
No. 20070420125 and No. 200801230. The author also gratefully
acknowledges the support of K. C. Wong Education Foundation, Hong
Kong}
%%\subjclass{}%
\thanks{E-mail:
xwchen$\symbol{64}$mail.ustc.edu.cn}
\maketitle
\date{}%
\dedicatory{}%
\commby{}%
\begin{center}
\end{center}

\begin{abstract}
We give a short proof to the following tilting theorem by Happel,
Reiten and Smal{\o} via an explicit construction: given two abelian
categories $\mathcal{A}$ and $\mathcal{B}$ such that $\mathcal{B}$
is tilted from $\mathcal{A}$, then $\mathcal{A}$ and $\mathcal{B}$
are derived equivalent.
\end{abstract}

\section{Introduction}
Let $\mathcal{A}$ be an abelian category. Recall that a
\emph{torsion pair} on $\mathcal{A}$ means a pair
$(\mathcal{T},\mathcal{F})$ of full subcategories  satisfying
\\
(T1).\quad ${\rm Hom}_\mathcal{A}(T, F)=0$ for all $T\in
\mathcal{T}$ and $F\in \mathcal{F}$; both subcategories
$\mathcal{T}$ and $\mathcal{F}$ are closed
under direct summands; \\
(T2). \quad  for each  object $X\in \mathcal{A}$, there is a short
exact sequence $0\longrightarrow T \longrightarrow X \longrightarrow
F \longrightarrow 0$ for some $T\in \mathcal{T}$ and $F\in
\mathcal{F}$.\par \vskip 5pt

In a torsion pair $(\mathcal{T}, \mathcal{F})$, it follows that the
subcategory $\mathcal{T}$ is closed under extensions and factor
objects; $\mathcal{F}$ is closed under extensions and sub objects.
The torsion pair $(\mathcal{T}, \mathcal{F})$ is called a
\emph{tilting torsion pair}, provided that each object in
$\mathcal{A}$ embeds into an object in $\mathcal{T}$. Dually the
torsion pair $(\mathcal{T}, \mathcal{F})$ is called a
\emph{cotilting torsion pair}, provided that each object in
$\mathcal{A}$ is a factor object of an object in $\mathcal{F}$
(\cite[Chapter I, section 3]{HRS}).

\par \vskip 5pt
Denote a complex in $\mathcal{A}$ by  $X^\bullet=(X^n, d_X^n)_{n\in
\mathbb{Z}}$ where $d_X^n: X^n\longrightarrow X^{n+1}$ is the
differential satisfying $d_X^{n+1}\circ d_X^n=0$; its shift
$X^\bullet[1]$ is a complex given by $(X^\bullet[1])^n=X^{n+1}$ and
$d_{X[1]}^n=-d_X^{n+1}$. Denote by $D(\mathcal{A})$ the (unbounded)
\emph{derived category} of $\mathcal{A}$, $D^{+}(\mathcal{A})$,
$D^{-}(\mathcal{A})$ and $D^b(\mathcal{A})$ the full subcategory
consisting of bounded-below, bounded-above and  bounded complexes,
respectively (\cite{V, Har}). We will always identify the abelian
category $\mathcal{A}$ as the full subcategory of $D(\mathcal{A})$
consisting of stalk complexes concentrated at degree zero
(\cite[p.40, Proposition 4.3]{Har}). Let $(\mathcal{T},
\mathcal{F})$ be a torsion pair on $\mathcal{A}$. Following
\cite[Chapter I, section 2]{HRS}, set $\mathcal{B}$ to be the full
subcategory of $D(\mathcal{A})$ consisting of complexes $X^\bullet$
satisfying $H^0(X^\bullet)\in \mathcal{T}$, $H^{-1}(X^\bullet)\in
\mathcal{F}$ and $H^i(X^\bullet)=0$ for $i\neq 0,1$. Note that
$\mathcal{T}\subseteq \mathcal{B}$ and $\mathcal{F}[1]\subseteq
\mathcal{B}$. By \cite[Chapter I, Proposition 2.1]{HRS} the category
$\mathcal{B}$ is the heart of certain $t$-structure on
$D(\mathcal{A})$ and thus by \cite{BBD} it is an abelian category
(also see \cite{GM2}); moreover the pair $(\mathcal{F}[1],
\mathcal{T})$ is a torsion pair on $\mathcal{B}$.
\par \vskip 5pt

 One might expect that the resulting new abelian category
$\mathcal{B}$ is derived equivalent to $\mathcal{A}$. But this is in
general false (by the example in \cite[p.16]{HRS}). However Happel,
Reiten and Smal{\o} show the following remarkable result
(\cite[Chapter I, Theorem 3.3]{HRS}).

\vskip10pt

\noindent {\bf Theorem} (Happel-Reiten-Smal{\o}) \emph{Let
$(\mathcal{T}, \mathcal{F})$ be a tilting torsion pair on
$\mathcal{A}$. Then we have a natural equivalence of triangulated
categories $D(\mathcal{B})\simeq D(\mathcal{A})$ which is compatible
with the inclusion of $\mathcal{B}$ into $D(\mathcal{A})$. Similar
results hold for $D^{*}(-)$ with $*\in \{+, -, b\}$.}

\vskip 10pt

 In the case of Theorem the category $\mathcal{B}$ is said to be \emph{tilted} from $\mathcal{A}$.
 Note that the original theorem only claims the
equivalence between the bounded derived categories and requires the
existence of enough projective or injective objects. The quoted
version is improved by Noohi (\cite[Theorem 7.6]{Noo}). We will give
a short proof of the theorem via an explicit construction of the
equivalence functor.

\section{The Proof of Theorem}

Throughout  $(\mathcal{T}, \mathcal{F})$ is a torsion pair on
$\mathcal{A}$ and $\mathcal{B}$ is the resulting abelian category.
We start with an easy observation.

\begin{lem}
Consider an complex $T^\bullet=(T^n,d_T^n)_{n\in \mathbb{Z}}$ with
terms in $\mathcal{T}$. Then it is exact in $\mathcal{A}$ if and
only if it is exact in $\mathcal{B}$.
\end{lem}

\noindent {\bf Proof.}\quad Assume that $T^\bullet$ is exact in
$\mathcal{A}$. Since $\mathcal{T}$ is closed under factor objects in
$\mathcal{A}$, the complex $T^\bullet$ splits into short exact
sequences $\xi^n: 0\longrightarrow T'^n
\stackrel{i^n}\longrightarrow T^n \stackrel{p^n}\longrightarrow
T'^{n+1} \longrightarrow 0$ with $T'^n\in \mathcal{T}$ and
$d_T^n=i^{n+1}\circ p^n$. Since $\mathcal{B}$ is a heart of certain
$t$-structure on $D(\mathcal{A})$, a sequence $0\longrightarrow
B^0\stackrel{f}\longrightarrow B^1 \stackrel{g}\longrightarrow
B^2\longrightarrow 0$ in $\mathcal{B}$ is short exact if and only if
there is a triangle $B^0\stackrel{f}\longrightarrow B^1
\stackrel{g}\longrightarrow B^2\longrightarrow B^0[1]$ in
$D(\mathcal{A})$ (\cite{BBD} and \cite[Chapter IV, \S 4]{GM2}). Note
that short exact sequences in $\mathcal{A}$ induces triangles in
$D(\mathcal{A})$ (\cite[p.62, Proposition 6.1]{Har}). Hence $\xi^n$
become short exact sequences in $\mathcal{B}$. Thus by splicing them
together we show that the complex $T^\bullet$ is exact in
$\mathcal{B}$. The ``if" part is proved similarly. \hfill
$\blacksquare$

\vskip 5pt

 The following result  is needed.

\begin{lem} {\rm (the ``only if" part of \cite[Chapter I, Proposition 3.2 \emph{(i)}]{HRS})}
Let $(\mathcal{T}, \mathcal{F})$ be a torsion pair on $\mathcal{A}$
and let $\mathcal{B}$ as before. If the torsion pair $(\mathcal{T},
\mathcal{F})$ is tilting, then the resulting torsion pair
$(\mathcal{F}[1], \mathcal{T})$ on $\mathcal{B}$ is cotilting.
\end{lem}

\vskip 5pt

\begin{rem}Note that the converse of the above lemma is also true as
stated in \cite[Chapter I, Proposition 3.2 \emph{(i)}]{HRS}. However
it seems to the author that a dual argument of this lemma is not
working. Instead, thanks to Theorem and then by combining
\cite[Chapter I, Proposition 3.4]{HRS} and the ``only if" part of
\cite[Chapter I, Proposition 3.2 \emph{(ii)}]{HRS} one deduces that
the converse holds (here one needs the fact that the equivalence in
Theorem is compatible with the inclusion $\mathcal{B}\hookrightarrow
D(\mathcal{A})$).
\end{rem}
\vskip 10pt

\noindent{\bf Proof of Theorem:}\quad  Denote by $K(\mathcal{A})$
the homotopy category of complexes in $\mathcal{A}$,
$K(\mathcal{T})$ (\emph{resp.} $K_{\rm ex}(\mathcal{A})$) its full
subcategory consisting of complexes in $\mathcal{T}$ (\emph{resp.}
exact complexes). The inclusion $K(\mathcal{T})\hookrightarrow
K(\mathcal{A})$ induces the following exact functor
$$F: K(\mathcal{T})/{K(\mathcal{T})\cap K_{\rm ex}(\mathcal{A})} \longrightarrow D(\mathcal{A}). $$
Since the torsion pair $(\mathcal{T}, \mathcal{F})$ is tilting and
$\mathcal{T}$ is closed under factor objects, we have for each $X\in
\mathcal{A}$ a short exact sequence $0\longrightarrow X
\longrightarrow T^0\longrightarrow T^1 \longrightarrow 0$ with
$T^i\in \mathcal{T}$. Note further that $\mathcal{T}$ is closed
under extensions, we infer that the conditions in \cite[p.42, Lemma
4.6 2)]{Har} are fulfilled, and thus for each  complex $X^\bullet$
in $K(\mathcal{A})$ there is a quasi-isomorphism $X^\bullet
\longrightarrow T^\bullet$ with $T^\bullet \in K(\mathcal{T})$. This
implies that the functor $F$ is dense and by \cite[p.283, 4-2
Th\'{e}or\`{e}me]{V} it is fully-faithful, that is, the functor $F$
is an equivalence of triangulated categories. By Lemma 2.2 we may
apply the dual argument to obtain a natural equivalence $$G:
K(\mathcal{T})/{K(\mathcal{T})\cap K_{\rm ex}(\mathcal{B})}
\longrightarrow D(\mathcal{B}).$$
 By Lemma 2.1 we have
$K(\mathcal{T})\cap K_{\rm ex}(\mathcal{A})=K(\mathcal{T})\cap
K_{\rm ex}(\mathcal{B})$. Hence
$FG^{-1}:D(\mathcal{B})\longrightarrow D(\mathcal{A})$ is the
required equivalence, where $G^{-1}$ denotes a quasi-inverse of $G$.
\par \vskip 5pt

To see other equivalences, let $*\in \{+,-,b\}$ and let $K^*(-)$
denote the corresponding homotopy categories. Note that in the
argument above, for a complex $X^\bullet\in K^*(\mathcal{A})$ we may
take a quasi-isomorphism $X^\bullet \longrightarrow T^\bullet$ with
$T^\bullet\in K^*(\mathcal{T})$ (for the case $*=+$, just consult
the proof in \cite[p.43, 1)]{Har};  for the case $*=-$, because
$\mathcal{T}$ is closed under factor objects one may replace
$T^\bullet$ by its good truncations; for the case $*=b$, consult the
proof in \cite[p.43, 1)]{Har} and note that since $\mathcal{T}$ is
closed under factor objects,  the argument therein is done within
finitely many steps, consequently the obtained complex $T^\bullet$
is bounded). Thus we construct the equivalences $F^*$ and $G^*$ as
above. This proves the corresponding equivalences between the
derived categories $D^*(-)$.
\par \vskip 5pt

 Finally we will show that the obtained equivalence  $FG^{-1}$ is compatible with the inclusion
$\mathcal{B}\hookrightarrow D(\mathcal{A})$. This is subtle. Given
an object $B\in \mathcal{B}$, since the torsion pair
$(\mathcal{F}[1], \mathcal{T})$ is cotilting, we have a short exact
sequence in $\mathcal{B}$, $\eta: 0\longrightarrow T^{-1}
\stackrel{d}\longrightarrow T^0\stackrel{g} \longrightarrow B
\longrightarrow 0$ with $T^i\in \mathcal{T}$, in other words, a
triangle $\xi: T^{-1} \stackrel{d}\longrightarrow
T^0\stackrel{g}\longrightarrow B \longrightarrow T^{-1}[1]$ in
$D(\mathcal{A})$. Then by construction $FG^{-1}(B)$ is isomorphic to
the complex $T^\bullet= \cdots \longrightarrow 0\longrightarrow
T^{-1}\stackrel{d}\longrightarrow T^0\longrightarrow
0\longrightarrow \cdots$. Note that the complex $T^\bullet$ is the
mapping cone of $d$ and thus form the triangle $\xi$ we obtain
$T^\bullet$ is isomorphic to $B$ (\cite[p.23, Propostion 1.1
c)]{Har}), in particular $T^\bullet \in \mathcal{B}$. Note the
following natural triangle $T^{-1}\stackrel{d}\longrightarrow T^0
\longrightarrow T^\bullet \longrightarrow T^{-1}[1]$ and thus a
short exact sequence $\gamma: 0\longrightarrow T^{-1}
\stackrel{d}\longrightarrow T^0\longrightarrow T^\bullet
\longrightarrow 0$ in $\mathcal{B}$. Comparing the short exact
sequences $\eta$ and $\gamma$ we obtain  a unique isomorphism
$\theta_B: B\simeq T^\bullet$ in $\mathcal{B}$. We claim that
$\theta$ is natural in $B$ and then we obtain  a natural isomorphism
between the inclusion functor $\mathcal{B}\hookrightarrow
D(\mathcal{A})$ and the composite $\mathcal{B}\hookrightarrow
D(\mathcal{B}) \stackrel{FG^{-1}}\longrightarrow D(\mathcal{A})$
(here we identify $T^\bullet$ with $FG^{-1}(B)$). \par \vskip 5pt

In fact, given a morphism $f: B\longrightarrow B'$ in $\mathcal{B}$,
choose an exact sequence $\eta': 0\longrightarrow T'^{-1}
\stackrel{d'}\longrightarrow T'^0\stackrel{g'} \longrightarrow B'
\longrightarrow 0$ with $T'^i\in \mathcal{T}$. Form the complex
$T'^\bullet$ and then obtain the short exact sequence $\gamma'$ and
the isomorphism $\theta_{B'}$ as above. Identify $G(T^\bullet)$ with
$B$, $G(T'^\bullet)$ with $B'$. Since the functor $G$ is
fully-faithful, we have a chain map $\phi^\bullet: T^\bullet
\longrightarrow T'^\bullet$ such that $G(\phi^\bullet)=f$. This
implies the following commutative exact diagram in $\mathcal{B}$

\[\xymatrix{
0\ar[r] & T^{-1} \ar[d]^{\phi^{-1}} \ar[r]^{d} & T^0 \ar[d]^{\phi^0}
\ar[r]^{g} & \ar[d]^{f} \ar[r] & 0\\
0 \ar[r] & T'^{-1} \ar[r]^{d'} & T'^0 \ar[r]^{g'} &B' \ar[r] & 0.
}\]
From this it is direct to see that $\theta_{B'}\circ
f=\phi^\bullet\circ \theta_B $ in $\mathcal{B}$ and thus in
$D(\mathcal{A})$. This finishes the proof.
 \hfill $\blacksquare$

\bibliography{}

\vskip 20pt

 {\footnotesize \noindent Xiao-Wu Chen, Department of
Mathematics, University of Science and Technology of
China, Hefei 230026, P. R. China \\
\emph{Current address}: Institut fuer Mathematik, Universitaet
Paderborn, 33095, Paderborn, Deutschland}

\end{document}